\documentclass[11pt]{amsart}
\usepackage{graphicx}
\usepackage{amsmath,amssymb}

\numberwithin{equation}{section}

\newtheorem{theorem}{Theorem}[section]
\newtheorem{corollary}[theorem]{Corollary}
\newtheorem{lemma}[theorem]{Lemma}

\newtheorem{remark}[theorem]{Remark}
\newtheorem{definition}[theorem]{Definition}
\begin{document}

\title[On Kaplansky's sixth conjecture]{On Kaplansky's sixth conjecture}
\thanks{This work was partially
supported by Natural Science Foundation of China (Grant No.
11201231), the China Postdoctoral Science Foundation (Grant No.
2012M511643), the Jiangsu Planned Projects for Postdoctoral Research
Funds (Grant No. 1102041C) and Agricultural Machinery Bureau
Foundation of Jiangsu Province (Grant No. GXZ11003).}

\address{College of Engineering, Nanjing Agricultural
University, Nanjing 210031, Jiangsu, China}

\author[L. Dai]{Li Dai}
\email[L. Dai]{daili1980@njau.edu.cn}
\author[J. Dong]{Jingcheng Dong}
\email[J. Dong]{dongjc@njau.edu.cn}



\begin{abstract}
About $39$ years ago, Kaplansky conjectured that the
dimension of a semisimple Hopf algebra over an algebraically closed
field of characteristic zero is divisible by the dimensions of its
simple modules. Although it still remains open, some partial answers
to this conjecture play an important role in classifying semisimple
Hopf algebras. This paper focuses on the recent development of
Kaplansky's sixth conjecture and its applications in classifying
semisimple Hopf algebras.
\end{abstract}

\maketitle

\section{Introduction}\label{sec1}
Let $H$ be a finite-dimensional semisimple Hopf algebra over an algebraically
closed field $k$. Kaplansky conjectured that the dimension of every simple
$H$-module divides the dimension of $H$. This is the sixth of a list
of ten conjectures posed by Kaplansky in his lecture notes ``Bialgebras"
\cite{Kaplansky}. Unfortunately, this conjecture is false
even for group algebras, which was already known at that time. For
example, let $G$ be the special linear group $SL(2,p)$ of $2\times
2$-matrices over a field with $p$ elements, where $p$ is an odd
prime, and let $kG$ be the group algebra of $G$ over an
algebraically closed field of characteristic $p$. Then $kG$ has
simple modules  whose dimensions do not divide the order of
$G$. See \cite[Example 17.17]{Curtis} for details. Of course the conjecture for group algebras holds true when the characteristic of the field $k$ is $0$ (thanks to a well-known result of Frobenius). We therefore believe that Kaplansky had ${\rm char}k=0$ in mind, but the published version missed it. So we can rewrite the conjecture as follows:

Let $H$ be a finite-dimensional semisimple Hopf algebra over an
algebraically closed field of characteristic $0$. Then the dimension
of a simple $H$-module divides the dimension of $H$.

We say that a semisimple Hopf algebra is of Frobenius type if it satisfies the conjecture above, in honour of Frobenius for his work.

The work toward solving Kaplansky's sixth conjecture can be roughly divided into two directions. The first direction is to consider semisimple Hopf algebras with low-dimensional simple modules. The pioneering work in this direction was made by Nichols and Richmond \cite{Nichols}. They proved, by analyzing the character algebra of a semisimple Hopf algebra, that if a semisimple Hopf algebra has a $2$-dimensional simple module then $2$ divides the dimension of the Hopf algebra. Their work has motivated great interest in this field which has produced many nice results. For example, using a similar method, Burciu \cite{Burciu}, Dong  and Dai \cite{dodai},  and Kashina et al
\cite{Kashina2} independently proved that if an odd-dimensional semisimple Hopf algebra has a $3$-dimensional simple module then $3$ divides the dimension of the Hopf algebra.

Their technique is also used to determine whether a low-dimensional semisimple Hopf algebra is of Frobenius type, since such a Hopf algebra often has a simple module of low dimension. For example, Natale \cite{Natale4} and Kashina \cite{Kashina3} independently proved that semisimple Hopf algebras of dimension less than $60$  are of Frobenius type.

The second direction is to think about semisimple Hopf algebras with particular properties. For example, Etingof and Gelaki proved that any quasitriangular semisimple Hopf algebra satisfies Kaplansky's sixth conjecture \cite{Etingof}. Some other results in this direction were made by Zhu \cite{ZhuS} for semisimple Hopf algebras whose characters are central in $H^*$, Zhu \cite{Zhu3} for semisimple Hopf algebras with a transitive module algebra, and Montgomery and Witherspoon \cite{Montgomery2} for semisolvable semisimple Hopf algebras.

This direction is also tightly related to the classification of semisimple Hopf algebras of a given dimension. As we will discuss in Section \ref{sec3}, semisimple Hopf algebras of dimension $p^n,pq,pq^2$ and $pqr$, where $p,q,r$ are distinct prime numbers and $n=1,2$ or $3$, have been completely classified. All these semisimple Hopf algebras are of Frobenius type. We will discuss in detail the recent work made by Etingof, Nikshych and Ostrik \cite{Etingof2} which covers all dimensions mentioned above.

Many examples show that a positive answer to Kaplansky's sixth
conjecture would be very helpful in classifying semisimple Hopf algebras.
For example, in the case that $H$ is a semisimple Hopf
algebra of dimension $pq$, where $p,q$ are distinct prime numbers, Gelaki and
Westreich \cite{GeWe} proved that if $H$ and $H^*$ are both of Frobenius type then $H$ is trivial; that is, it is either a group algebra or a dual group algebra. In a subsequent paper \cite{EtingofG}, Etingof and Gelaki proved that $H$ and $H^*$ are of Frobenius type, and hence completed the classification of $H$. This result was also obtained by Sommerh\"{a}user \cite{Somm} by different methods. Another example is taken from Natale's work. Let $H$ be a semisimple Hopf algebra of dimension $pq^2$, where $p,q$ are distinct prime numbers. In
\cite{Natale3}, Natale completed the classification of $H$ by assuming that $H$ and $H^*$ are both of Frobenius type. Some other applications of Kaplansky's sixth
conjecture may be found in the authors' recent work \cite{dong}, \cite{dong1}, \cite{dodai}, \cite{dong2}.

There are three nice reviews related to our subject \cite{Burciu2}, \cite{Natale5}, \cite{Sommer}. In \cite{Sommer}, Sommerh\"{a}user reviewed all of Kaplansky's ten conjectures. In \cite[Section 1]{Burciu2}, Burciu reviewed the results on Kaplansky's sixth conjecture obtained until then, and mainly focused on the development of semisimple Hopf algebras. In \cite[Section 6]{Natale5}, Natale gave a brief review on Kaplansky's sixth conjecture and mainly paid attention to the representations of semisimple Hopf algebras.

In the fusion category setting, there is a similar question: Is every fusion category of Frobenius type? Here, a fusion category $\mathcal{C}$ is of Frobenius
type if for every simple object $X$ of $\mathcal{C}$, the Frobenius-Perron dimension ${\rm FPdim} X$ of $X$ divides the Frobenius-Perron dimension ${\rm FPdim} \mathcal{C}$ of $\mathcal{C}$; that is, the ratio ${\rm FPdim} \mathcal{C}/{\rm FPdim} X$ is an algebraic integer. We will review in Subsection \ref{sec2-2} and Subsection \ref{sec3-5} some recent developments in this direction.

In this article we shall review results and approaches so far in
the study of Kaplansky's sixth conjecture, as well as its applications in
classifying semisimple Hopf algebras. In the last part of  this
article we shall also present our point of view on solving this
conjecture.

Throughout, we will work over an algebraically closed field $k$ of
characteristic $0$. Our references for the theory of Hopf algebras
are \cite{Montgomery} or \cite{Sweedler2}.

\section{Low-dimensional simple modules and semisimple Hopf algebras}\label{sec2}
\subsection{Semisimple Hopf algebras}
A Hopf algebra is called semisimple (respectively, cosemisimple) if
it is semisimple as an algebra (respectively, if it is cosemisimple
as a coalgebra). A semisimple Hopf algebra is automatically
finite-dimensional by \cite[Corollary 2.7]{Sweedler}. By a result of
Larson and Radford \cite{Larson}, \cite{Larson2}, a finite-dimensional Hopf
algebra is semisimple if and only if it is cosemisimple.

Let $H$ be a semisimple Hopf algebra and let $V$ be an $H$-module.
The character of $V$ is the element $\chi_V\in H^*$ defined by
$\langle\chi_V,h\rangle={\rm Tr}_V(h)$ for all $h\in H$. The degree
of $\chi_V$ is defined to be the integer ${\rm
deg}\chi_V=\chi_V(1)={\rm dim}V$.

By a result of Zhu \cite{Zhu}, the irreducible characters of $H$,
namely, the characters of the simple $H$-modules, span a semisimple
subalgebra $R(H)$ of $H^*$, which is called the character algebra of
$H$. The antipode $S$ induces an anti-algebra involution $*: R(H)\to
R(H)$, given by $\chi\mapsto\chi^*:=S(\chi)$. We call $\chi^*$ the
dual of $\chi$. Let ${\rm Irr}(H)$ denote the set of non-isomorphic irreducible characters of $H$. Then ${\rm Irr}(H)$ is a $k$-basis of $R(H)$.

Pioneers in solving Kaplansky's sixth conjecture, Nichols
and Richmond began their work by considering semisimple Hopf
algebras with simple modules of dimension $2$. They proved
\cite[Theorem 11]{Nichols}:

\begin{theorem}\label{thm1}
If a semisimple Hopf algebra $H$ has a simple module of dimension $2$
then the dimension of the semisimple Hopf algebra is even.
\end{theorem}

Besides the importance of the result itself, the technique used in
\cite{Nichols} is also important. To prove their main result,
Nichols and Richmond analyzed  the possible decomposition of
$\chi\chi^*$, where $\chi$ is an irreducible character of degree
$2$, and tried to look for standard subalgebras of
$R(H)$. Here, a standard subalgebra of
$R(H)$ is a subalgebra of $R(H)$ which is spanned by a subset of the basis ${\rm Irr}(H)$. Their main result then follows from the following theorem
\cite[Theorem 6]{Nichols}:

\begin{theorem}\label{thm1-1}
There is a bijection between standard subalgebras of $R(H)$ and
quotient Hopf algebras of $H$.
\end{theorem}

They finally proved that $H$ admits certain quotient
Hopf algebras of dimension $2$, $12$, $24$ or $60$. Therefore, the
dimension of $H$ is even, in the light of the main
theorem in \cite{Nichols2}.

About six years later, Kashina, Sommerh\"{a}user and  Zhu generalized
the above result. They proved \cite[Theorem 4.1, Theorem
5.1]{Kashina}:

\begin{theorem}\label{thm2}
Let $H$ be a semisimple Hopf algebra. Then

(1)\, If $H$ has a non-trivial self-dual simple module, then the dimension of $H$ is even.

(2)\, If $H$ has a simple module of even dimension, then the dimension of $H$ is even.
\end{theorem}

The proof of the first part is heavily dependent on the Frobenius-Schur theorem for Hopf algebras \cite{Linchenko} and the result on the exponent of a semisimple Hopf algebra \cite[Proposition 2.2]{Kashina}, while the second part is based on the first part and an analysis of the decomposition of the product of an irreducible character of even degree and its dual.

Theorem \ref{thm2} is also very important in studying Kaplansky's sixth conjecture. We shall discuss it at the end of this subsection.

Recently, Bichon and Natale gave a more precise description of the work of Nichols and Richmond \cite[Theorem 1.1]{BiNa2011}. They proved:

\begin{theorem}
Let $H$ be a cosemisimple Hopf algebra. Suppose that $H$ has an irreducible cocharacter $\chi$ of degree $2$ and $C$ is the simple subcoalgebra containing $\chi$. Then the
subalgebra $B = k[CS(C)]$ is a commutative Hopf subalgebra of $H$ isomorphic
to $k^{G}$, where $G$ is a non-cyclic finite subgroup of ${\rm PSL}_2(k)$ of even order.

More specifically, let $G[\chi] \subseteq G(H)$ be the stabilizer of $\chi$ under the left multiplication by $G(H)$, then the order of $G[\chi]$ divides $4$, and the following hold:

(1) If the order of $G[\chi]$ is 4, then $B\cong k^{\mathbb{Z}_2\times\mathbb{Z}_2}$.

(2) If the order of $G[\chi]$ is 2, then $B\cong k^{D_n}$, where $n\geq3$.

(3) If the order of $G[\chi]$ is 1, then $B\cong k^{\mathbb{A}_4}, k^{\mathbb{S}_4}$, or $k^{\mathbb{A}_5}$.
\end{theorem}

 They also studied the special case when the $2$-dimensional simple comodule is faithful, and more interesting results were obtained. They then applied their results to the classification of semisimple Hopf algebras of dimension $60$ and of semisimple Hopf algebras such that the dimensions of its simple comodules are at most $2$.

\medbreak
Motivated by the work of Nichols and Richmond, many algebraists
intend to consider semisimple Hopf algebras with
a simple module of dimension $3$.  Unfortunately, only partial
answers have been obtained. For example, Burciu
\cite{Burciu}, Dong  and Dai \cite{dodai} and Kashina et al
\cite{Kashina2} independently proved:

\begin{theorem}\label{thm3}
If a semisimple Hopf algebra is of odd dimension and has a simple
module of dimension $3$, then the dimension of the Hopf algebra is
divisible by $3$.
\end{theorem}

All these three results are mainly based on a similar treatment as
in \cite{Nichols}, while the last two articles also adopt the main
result in \cite{Kashina} which greatly simplifies the proof.

Besides the applications above, the technique used by Nichols and
Richmond is also very useful in determining whether a semisimple
Hopf algebra of low dimension satisfies Kaplansky's sixth conjecture,
because a semisimple Hopf algebra of low dimension often admits a
simple module of low dimension, such as $2$ or $3$.

Following the technique in \cite{Nichols}, Natale \cite{Natale4} and
Kashina \cite{Kashina3} independently proved that semisimple Hopf
algebras of dimension less than $60$  satisfy Kaplansky's sixth
conjecture. Moreover, Natale took a further step to prove that all
these Hopf algebras are either upper or lower semisolvable up to a
cocycle twist. The notion of semisolvability will be given in the next section. Therefore, Natale completes the classification of all
these Hopf algebras, to some degree.

Now we illustrate why Theorem \ref{thm2} is important in studying Kaplansky's sixth conjecture. Let $H$ be a semisimple Hopf algebra over $k$. As an algebra, $H$ is isomorphic to a direct product of full matrix algebras
$$H\cong k^{(n_1)}\times \prod_{i=2}^{s}M_{d_i}(k)^{(n_i)},$$

where $n_1=|G(H^*)|$. In this case, $H$ is called of type $(d_1,n_1;\cdots;d_s,n_s)$ as an
algebra, where $d_1=1$. Obviously, $H$ is of type $(d_1,n_1;\cdots;d_s,n_s)$ as an algebra if and only if $H$ has $n_1$ non-isomorphic irreducible characters of degree $d_1$, $n_2$
non-isomorphic irreducible characters of degree $d_2$, and so on.

Suppose that the dimension of $H$ is odd and $H$ is of type
$(d_1,n_1;\cdots;d_s,n_s)$ as an algebra. Then part (2) of Theorem
\ref{thm2} clearly shows that $d_i$ is odd,  and part (1) of Theorem
\ref{thm2} shows that $n_i$ is even for all $2\leq i\leq s$. Indeed,
if there exists $i\in \{2,\cdots,s\}$ such that $n_i$ is odd, then
there is at least one irreducible character of degree $d_i$ which is self-dual. This contradicts part (1) of Theorem \ref{thm2}. Therefore, Theorem
\ref{thm2} is  quite useful in excluding potential type
$(d_1,n_1;\cdots;d_s,n_s)$ for a semisimple Hopf algebra. Using this
observation, together with other techniques, Dong and
Dai \cite{dodai} further extended the results of Natale \cite{Natale4}
and Kashina \cite{Kashina3}. That is, they proved that
odd-dimensional semisimple Hopf algebras of dimension less than
$600$ satisfy Kaplansky's sixth conjecture.

We should remark that the technique used in \cite{Nichols} does not work well for simple modules of higher dimension.

\subsection{Results from fusion categories}\label{sec2-2}
A fusion category over $k$ is a $k$-linear semisimple rigid tensor
category with finitely many isomorphism classes of simple objects,
finite-dimensional spaces of morphisms, and simple unit object. Let $H$ be a semisimple Hopf algebra over $k$. Then the category ${\rm Rep}H$ of its finite-dimensional representations is a fusion
category.

Recall that a fusion category is said to be weakly integral if its Frobenius-Perron dimension is an integer. A fusion category is said to be integral if the Frobenius-Perron dimension of every simple object is an integer.

\begin{theorem}
Let $\mathcal{C}$ be an integral fusion category. Suppose that the Frobenius-Perron dimensions of its simple objects are $1,2$ or $3$. Then $\mathcal{C}$ is of Frobenius type.
\end{theorem}

The theorem above is the main result of \cite{dong2013integral}. Its proof relies on an analogue of Theorem \ref{thm1-1} in the fusion category setting. The proof of Theorem \ref{thm1-1} only makes use of properties of the Grothendieck ring of a semisimple Hopf algebra. Therefore its proof also works  \textit{mutatis mutandis} in the fusion category setting.

In \cite{Etingof2}, Etingof, Nikshych and Ostrik proved that any weakly integral fusion category of Frobenius-Perron dimension less than $84$ is of Frobenius type. The following theorem \cite[Theorem 1.1]{dong2012frobenius} extends this result.

\begin{theorem}\label{thm1-2}
Let $\mathcal{C}$ be a weakly integral fusion category of Frobenius-Perron dimension less than $120$. Then $\mathcal{C}$ is of Frobenius type. Furthermore,  if ${\rm FPdim} \mathcal{C} > 1$ and $\mathcal{C} \ncong{\rm Rep}\mathbb A_5$, then $\mathcal{C}$ has nontrivial invertible objects.
\end{theorem}

A fusion category is called simple if it has no nontrivial proper fusion
subcategories \cite{Etingof2}.  As a consequence of Theorem \ref{thm1-2}, if
${\rm FPdim} \mathcal{C} \leq 119$ and ${\rm FPdim} \mathcal{C} \neq 60$ or $p$, where $p$ is
a prime number, then $\mathcal{C}$ is not simple as a fusion
category. Combined with the results of the paper \cite{Etingof2}, the
theorem above implies that the only weakly integral simple fusion
categories of Frobenius-Perron dimension $\leq 119$ are the
categories ${\rm Rep} \mathbb A_5$ of finite-dimensional
representations of the alternating group $\mathbb A_5$ and the
pointed fusion categories $\mathcal{C}(\mathbb Z_p, \omega)$ of
finite-dimensional $\mathbb Z_p$-graded vector spaces, where $p$
is a prime number, with associativity constraint determined by a
$3$-cocycle $\omega \in H^3(\mathbb Z_p, k^*)$.

\section{Semisimple Hopf algebras that satisfy Kaplansky's sixth conjecture}\label{sec3}
\subsection{Semisimple Hopf algebras whose characters are central in $H^*$}
Before the work of Nichols and Richmond, Zhu had already proven the following result on Kaplansky's sixth conjecture in \cite{ZhuS}.

\begin{theorem}\label{thm4}
Let $H$ be a  semisimple Hopf algebra. If $R(H)$ is central in $H^*$
then $H$ satisfies Kaplansky's sixth conjecture.
\end{theorem}

We refer to \cite{Lorenz} for an alternate proof of this theorem.
Although Zhu's result is interesting, except
for dual group algebras, we do not yet know which Hopf algebras
satisfy the assumptions.

\subsection{Quasitriangular semisimple Hopf algebras}

Let $H$ be a semisimple Hopf algebra. We define two actions of $H^*$
on $H$ as
$$f\rightharpoonup h=\sum f(h_2)h_1\mbox{\quad and\quad}
h\leftharpoonup f=\sum f(h_1)h_2,\mbox{\, for all \,}f\in H^*, h\in
H.$$

The Drinfeld double $D(H)$ of $H$ has $H^{*cop}\otimes H$ as
its underlying vector space with multiplication in $D(H)$ given
by
$$(g\otimes h)(f\otimes l)=\sum g(h_1\rightharpoonup f
\leftharpoonup S^{-1}(h_3))\otimes h_2l.$$
$D(H)$ has the coalgebra structure of the usual tensor product of
coalgebras. It follows from \cite{Montgomery} that $D(H)$ is also
semisimple. Etingof and Gelaki proved \cite[Theorem 1.4]{Etingof}:

\begin{theorem}\label{thm5}
If $H$ is a semisimple Hopf algebra and $V$ is a simple
$D(H)$-module, then the dimension of $V$ divides the dimension of
$H$.
\end{theorem}

This is a nice generalization of Zhu's work \cite{Zhu2} which states
that the dimensions of the simple $D(H)$-submodules of $H$ divide
the dimension $H$.

The proof of Theorem \ref{thm5} uses the Verlinde formula from
modular categories. A modular category is a fusion category with
nondegenerate $S$-matrix.

If $H$ is a quasitriangular semisimple Hopf algebra then the
universal $R$-matrix provides a surjective Hopf algebra map from
$D(H)$ to $H$. Therefore, every simple $H$-module is also a simple
$D(H)$-module via pull back, and hence the following result
\cite{Etingof} follows from Theorem \ref{thm5}:

\begin{theorem}\label{thm6}
If $H$ is a quasitriangular semisimple Hopf algebra then it
satisfies Kaplansky's sixth conjecture.
\end{theorem}

Two alternate proofs of  Theorem \ref{thm5} were later offered in
\cite{Sch} and \cite{Tsang} which both used the Class
Equation \cite{Lorenz2}, \cite{Zhu}. In addition, Schneider's article \cite{Sch} generalizes Theorem \ref{thm5} to factorizable Hopf
algebras.

\begin{theorem}[Class Equation]\label{thm6-1}
Let $H$ be a semisimple Hopf algebra and $R(H)$ its character algebra. For every primitive idempotent $e$ of $R(H)$, ${\rm dim}eH^*$ divides ${\rm dim}H$. If $e_1,\cdots,e_n$ are the primitive idempotents of $R(H)$, then
$${\rm dim}H=1+\sum_{i=2}^n{\rm dim}e_iH^*,$$
where $e_1$ is the normalized integral in $H^*$.
\end{theorem}

If $H=kG$ is a group algebra then the elements in $R(H)$ are constant on conjugacy classes $\mathcal{C}_1,\cdots,\mathcal{C}_n$. Let $G=\{g_1,\cdots,g_s\}$ and $\{p_{g_1}\cdots,p_{g_s}\}$ be the corresponding dual basis. Then $e_i=\sum_{g\in\mathcal{C}_i}p_g$, which implies that the size of every conjugacy class divides the order of $G$. So Theorem \ref{thm6-1} is the generalization of the usual Class Equation for finite groups.

If $H=k^G$ is a dual group algebra then $R(H)=kG$. Hence, $e_1,\cdots,e_n$ are the primitive idempotents of $kG$, and ${\rm dim}e_iH^*={\rm dim}e_ikG$ is the dimension of simple module associated to $e_i$. This is a well known result due to Frobenius which was mentioned in the Introduction.

The Class Equation is also used in the classification of semisimple Hopf algebras. We will elaborate the work of Zhu \cite{Zhu} and Masuoka \cite{Masuoka1}. Zhu proved the following theorem \cite{Zhu} which solves a conjecture of Kaplansky \cite{Kaplansky}. Similar ideas were used by Kac to get an analogous result in the setting of $C^*$-algebras \cite{Kac}.

\begin{theorem}
A Hopf algebra of prime dimension is necessarily semisimple and isomorphic to the group algebra $k(\mathbb{Z}/p\mathbb{Z})$, where $p$ is a prime number.
\end{theorem}

Masuoka later proved the following theorem \cite[Theorem 2]{Masuoka1} which was used to prove that a semisimple Hopf algebra of dimension $p^2$ is isomorphic
to the group algebra $k(\mathbb{Z}/p^2\mathbb{Z})$ or $k(\mathbb{Z}/p\mathbb{Z})^2$, where $p$ is a prime number.

\begin{theorem}
Suppose that the dimension of a semisimple Hopf algebra $H$ is $p^n$, where $p$ is a prime number and $n$ is a positive integer. Then $H$ has  a non-trivial central group-like element.
\end{theorem}

\medbreak
By the result of \cite{majid}, Theorem \ref{thm5} means that the
dimensions of simple Yetter-Drinfeld  $H$-modules divide the
dimension of $H$.

Note that a Yetter-Drinfeld submodule $M\subseteq H$ is exactly
a left coideal $M$ of $H$ such that $h_1M S(h_2)\subseteq M$ for
all $h\in H$. So, a $1$-dimensional Yetter-Drinfeld submodule of $H$
is exactly the span of a central group-like element of $H$. This
observation together with Theorem \ref{thm5} can be used to
determine the existence of normal Hopf subalgebras (a Hopf
subalgebra $A\subseteq H$ is called normal if $h_1AS(h_2)\subseteq
A$ and $S(h_1)Ah_2\subseteq A$ for all $h\in H$) as follows:

Let $\pi:H\to B$ be a Hopf algebra map and
$$H^{co\pi}=\{h\in H:(id\otimes \pi)\Delta(h)=h\otimes 1\}$$ be the
coinvariant subspace of $H$. Then $H^{co\pi}$ is a left coideal
subalgebra of $H$. Moreover, $H^{co\pi}$ is stable under the left
adjoint action of $H$ \cite{Sch2}. It follows that $H^{co\pi}$ is a
Yetter-Drinfeld submodule of $H$. Therefore, $H^{co\pi}$ is a direct
sum of simple Yetter-Drinfeld submodules of $H$ and the dimension of
every such simple module divides the dimension of $H$. By analyzing the
possible decompositions of $H^{co\pi}$ into simple Yetter-Drinfeld
submodules of $H$, we can determine whether $H$ contains central
group-like elements. This technique has been used in
\cite{dong}, \cite{dong1}, \cite{dodai}, \cite{dong2}, \cite{Natale4}.

\subsection{Semisolvable semisimple Hopf algebras}\label{sec4}

Let $H$ be a Hopf algebra, and let $A$ be an algebra. Suppose that $\sigma:H\otimes H\to A$ is a convolution-invertible $k$-linear map and
$\rightharpoonup:H\otimes A\to A$ is a $k$-linear map. Suppose further that, for every $h,l,m\in H,a,b\in A$, they satisfy:

(1)\, $h\rightharpoonup(l\rightharpoonup a)=\sum\sigma(h_1,l_1)(h_2l_2\rightharpoonup a)\sigma^{-1}(h_3,l_3)$;

(2)\, $h\rightharpoonup ab=\sum(h_1\rightharpoonup a)(h_2\to b),\,
h\rightharpoonup 1=\varepsilon(h)1$,\, $1\rightharpoonup a=a$;

(3)\,$\sum(h_1\rightharpoonup
\sigma(l_1,m_1))\sigma(h_2,l_2m_2)=\sum\sigma(h_1,l_1)\sigma(h_2l_2,m)$;

(4)\, $\sigma(h,1)=\varepsilon(h)1=\sigma(1,h)$.

Then the crossed product algebra $A\#_{\sigma}H$ is the
vector space $A\otimes H$ together with unit $1\otimes 1$ and the
multiplication
$$(a\#_{\sigma}h)(b\#_{\sigma} l)=a(h_1\rightharpoonup b)\sigma(h_2,l_1)\#_{\sigma}h_3l_2.$$

The notions of upper and lower semisolvability for finite-dimensional Hopf algebras were introduced in
\cite{Montgomery2}, as generalizations of the notion of solvability
for finite groups. By definition, $H$ is called lower semisolvable
if there exists a chain of Hopf subalgebras
$$H_{n+1} = k\subseteq H_{n}\subseteq\cdots \subseteq H_1 = H$$
such that $H_{i+1}$ is a
normal Hopf subalgebra of $H_i$, for all $i$, and all quotients
$H_{i}/H_{i}H^+_{i+1}$ are commutative or cocommutative. Dually, $H$ is called upper
semisolvable if there exists a chain of quotient Hopf algebras
$$H_{(0)} =
H\xrightarrow{\pi_1}H_{(1)}\xrightarrow{\pi_2}\cdots\xrightarrow{\pi_n}H_{(n)}
= k$$
such that $H_{(i-1)}^{co\pi_{i}}=\{h\in H_{(i-1)}|(id\otimes
\pi_i)\Delta(h)=h\otimes 1\}$ is a normal Hopf subalgebra of
$H_{(i-1)}$, and all $H_{(i-1)}^{co\pi_i}$ are commutative or cocommutative.

The following conjecture can be viewed as a generalization of
Kaplansky's sixth conjecture. When $H$ is lower or upper semisolvable it
was proved by Montgomery and Witherspoon \cite[Theorem
3.4]{Montgomery2}.

\textbf{Conjecture} If $A$ is a finite-dimensional semisimple
algebra of Frobenius type and $H$ is a semisimple Hopf algebra then
the crossed product $A\#_{\sigma}H$ is of Frobenius type.

The work of Montgomery and Witherspoon is very useful in determining
whether a non-simple semisimple Hopf algebra satisfies Kaplansky's sixth
conjecture. Note that a Hopf algebra is called simple if it does not
contain proper normal Hopf subalgebras. Indeed, let $A$ be a proper
normal Hopf subalgebra of $H$. Then by the result in \cite{Sch2},
$H\cong A\#_{\sigma}(H/HA^+)$ is a crossed product for some
$\sigma$. Therefore, if $A$ is of Frobenius type and $H/HA^+$ is
lower or upper semisolvable then $H$ satisfies Kaplansky's sixth
conjecture. Using this observation, Montgomery and Witherspoon
proved \cite[Theorem 3.5, Corollary 3.6]{Montgomery2}:

\begin{theorem}\label{thm7}
Let $H$ be a semisimple Hopf algebra of dimension $p^n$, where $p$
is a prime number and $n$ is an integer. Then $H$ is upper and lower
semisolvable and therefore satisfies Kaplansky's sixth conjecture.
\end{theorem}

In fact, when $n=1,2,3$ Theorem \ref{thm7} has been obtained in
\cite{Masuoka1}, \cite{Masuoka2}, \cite{Zhu} as a by-product of the classification of
semisimple Hopf algebras.

\subsection{Semisimple Hopf algebras with a transitive module algebra}
Let $H$ be a Hopf algebra. A module algebra of $H$ is an associative algebra $A$ on which $H$ acts via $h\cdot 1=\varepsilon(h)1$ and $h\cdot (ab)=\sum (h_1\cdot a)(h_2\cdot b)$, where $h\in H$ and $a,b\in A$.

$I\subseteq A$ is called a module ideal if $I$ is a two-sided ideal and $I$ is an $H$-submodule
of $A$. A module algebra $A$ of $H$ is called transitive if it satisfies the following two conditions:

(1)\,$A^H=\{a\in A| h\cdot a=\varepsilon(h)a,\, \mbox{for all\,\,} h\in H \}=k$;

(2)\,$A$ has no proper module ideals.

\medbreak
In \cite{YaZh}, \cite{Zhu3}, Yan and Zhu tried to solve Kaplansky's sixth conjecture by considering semisimple Hopf algebras with a transitive module algebra. In fact, every semisimple Hopf algebra has this property. For example, let $V$ be a simple $H$-module, and let $\alpha : H\to {\rm End}_k(V)$ be the corresponding algebra morphism. Considering the conjugation action of $H$ on ${\rm End}_k(V)$: $h\cdot f=\sum\alpha(h_1)f\alpha(S(h_2))$, ${\rm End}_k(V)$ becomes an $H$-module algebra. The simplicity of $V$ shows that ${\rm End}_k(V)$ is transitive. They proved \cite{Zhu3}:
\begin{theorem}\label{thm7-1}
Let $H$ be a semisimple Hopf algebra. If $A$ is
a transitive $H$-module algebra and $V$ is a simple
$A$-module, then ${\rm dim}A$ divides $({\rm dim}V)^2{\rm dim}H$.
\end{theorem}
Although the theorem above can not solve the conjecture, it is very close to that point. You may agree with this point of view by looking at the following conjecture \cite{Zhu3}:

\medbreak
\textbf{Conjecture}\quad  Let $H$ be a semisimple Hopf algebra and $A$ be a transitive
module algebra of $H$. Then for each simple $A$-module $V$, ${\rm dim}A$ divides ${\rm dim}V{\rm dim}H$.

\medbreak
In fact, this conjecture truly implies Kaplansky's sixth conjecture: Since the simple $H$-module $V$ is the unique simple ${\rm End}_k(V)$-module, the conjecture above means that ${\rm dim}V$ divides ${\rm dim}H$.

\subsection{Weakly group-theoretical semisimple Hopf algebras}\label{sec3-5}
Let $G$ be a finite group, and let $\mathcal{C}$ be a fusion category. A $G$-grading of $\mathcal{C}$  is a direct sum of full abelian subcategories  $\mathcal{C} = \oplus_{g \in G} \mathcal{C}_g$, where the tensor product of $\mathcal{C}$ maps $\mathcal{C}_g \times \mathcal{C}_h \to \mathcal{C}_{gh}$ and $(\mathcal{C}_g)^* = \mathcal{C}_{g^{-1}}$. The grading is called faithful if $\mathcal{C}_g\neq 0$, for all $g\in G$. In this case, $\mathcal{C}$ is called a $G$-extension of $\mathcal{C}_e$, where $\mathcal{C}_e$ is the neutral component of $\mathcal{C}$.

A fusion category $\mathcal{C}$ is said to be (cyclically) nilpotent if there is a sequence of
fusion categories $\mathcal{C}_0 ={\rm Vec}, \mathcal{C}_1, \cdots, \mathcal{C}_n =\mathcal{C}$ and a sequence of finite (cyclic) groups $G_1,\cdots, G_n$ such that $\mathcal{C}_i$ is obtained from $\mathcal{C}_{i-1}$ by a $G_i$-extension.

Let $G$ be a finite group and $\mathcal{C}$ be a fusion category. Let $\underline{G}$ denote the monoidal category whose objects are elements of $G$, morphisms are identities and the tensor product is given by the multiplication in $G$. Let ${\rm\underline{Aut}}_{\otimes}\mathcal{C}$ denote the monoidal category whose objects are tensor autoequivalences of $\mathcal{C}$, morphisms are isomorphisms of tensor functors and the tensor product is given by the composition of functors.

An action of $G$ on $\mathcal{C}$ is a monoidal functor
$$T:\underline{G}\to {\rm\underline{Aut}}_{\otimes}\mathcal{C},\quad g\mapsto T_g$$
with the isomorphism $f^V_{g,h}: T_g(T_h(V))\cong T_{gh}(V)$, for every $V$ in $\mathcal{C}$.

Let $\mathcal{C}$ be a fusion category with an action of $G$. Then the fusion category $\mathcal{C}^G$, called the $G$-equivariantization of $\mathcal{C}$, is defined as follows:

(1)\quad An object in $\mathcal{C}^G$ is a pair $(V, (u^V_g)_{g\in G})$, where $V$ is an object of $\mathcal{C}$ and $u^V_g: T_g(V)\to V$
is an isomorphism such that,
$$u^V_gT_g(u^V_h)= u^V_{gh}f^V_{g,h},\quad \mbox{for all\quad} g,h\in G.$$

(2)\quad A morphism $\phi: (U,u_g^U)\to (V,u_g^V)$ in $\mathcal{C}^G$ is a morphism $\phi: U\to V$ in $\mathcal{C}$ such that $\phi u_g^U=u_g^VT_g(\phi)$, for all $g\in G$.

(3)\quad The tensor product in $\mathcal{C}^G$ is  defined as $(U,u_g^U)\otimes (V,u_g^V)=(U\otimes V, (u_g^U\otimes u_g^V)j_g|_{U,V})$, where $j_g|_{U,V}:T_g(U\otimes V)\to T_g(U)\otimes T_g(V)$ is the isomorphism giving the monoidal structure of $T_g$.

Let $\mathcal{C},\mathcal{D}$ be fusion categories, and
$\mathcal{M}$ be an indecomposable left $\mathcal{C}$-module
category. Then $\mathcal{C}$ and $\mathcal{D}$ are Morita equivalent
if $\mathcal{D}$ is equivalent to $\mathcal{C}^*_{\mathcal{M}}$ which is the category of
$\mathcal{C}$-module endofunctors of $\mathcal{M}$.

As an analogue of the classical approach for algebras, we use
Morita equivalence to classify fusion categories.

A fusion category is called pointed if all of its simple objects are
invertible. A fusion category is called group-theoretical if
it is Morita equivalent to a pointed fusion category. A weakly group-theoretical fusion category is a fusion category which is Morita equivalent to a nilpotent fusion category.

\begin{definition}\label{def1}
A fusion category $\mathcal{C}$ is called solvable if it satisfies one of the following equivalent conditions:

(1)\quad $\mathcal{C}$ is Morita equivalent to a cyclically nilpotent fusion category;

(2)\quad There is a sequence of fusion categories $\mathcal{C}_0 ={\rm Vec}, \mathcal{C}_1, \cdots, \mathcal{C}_n =\mathcal{C}$ and a sequence of cyclic groups $G_1,\cdots, G_n$ of prime order such that $\mathcal{C}_i$ is obtained from $\mathcal{C}_{i-1}$ either by a $G_i$-equivariantization or by a $G_i$-extension.
\end{definition}

A semisimple Hopf algebra is called weakly group-theoretical
(respectively, solvable) if ${\rm Rep}H$ is weakly group-theoretical (respectively, solvable). We refer the reader to \cite{Etingof1}, \cite{Etingof2} for
further definitions and results about fusion categories.

We remark that solvability for semisimple Hopf algebras can also be
viewed as a generalization of the notion of solvability for finite
groups. But the interrelations between solvability and
semisolvability for semisimple Hopf algebras are still not clear.
The reader can find an explanation in \cite[Remark 4.6]{Etingof2}.

A fusion category $\mathcal{C}$ has the strong Frobenius
property if for every indecomposable $\mathcal{C}$-module category
$\mathcal{M}$ and any simple object $X$ in $\mathcal{M}$ the number
$\frac{{\rm FPdim}(\mathcal{C})}{{\rm FPdim}(X)}$ is an algebraic integer, where
the Frobenius-Perron dimension of $\mathcal{M}$ is normalized in
such a way that ${\rm FPdim}(\mathcal{M})= {\rm FPdim}(\mathcal{C})$. The strong
Frobenius property of a fusion category is a strong form of
Kaplansky's sixth conjecture. To see this, it suffices to take
$\mathcal{M}=\mathcal{C}$ the category of finite-dimensional
representations of a semisimple Hopf algebra.

Etingof, Nikshych and Ostrik proved the following theorem \cite[Theorem 1.5]{Etingof2}:

\begin{theorem}\label{thm8}
Any weakly group-theoretical fusion category has the strong
Frobenius property.
\end{theorem}

From the definitions above, we know that the class of weakly group-theoretical
fusion categories covers the classes of solvable and
group-theoretical fusion categories. Moreover, the class of weakly
group-theoretical fusion categories actually covers all known fusion
categories which are weakly integral.

The following two theorems \cite[Theorem 1.6, Theorem 9.2]{Etingof2}
are more concrete, and the first one is an analogue of Burnside's
theorem for fusion categories (Compare the question on
semisolvability for semisimple Hopf algebras posed by
Montgomery (2000), see \cite[Question 4.17]{Andr}).

\begin{theorem}\label{thm9}
(1)\, Any integral fusion category of Frobenius-Perron dimension $p^aq^b$ is
solvable, where $p,q$ are prime numbers and $a,b$ are non-negative
integers.

(2)\, Any integral fusion category of Frobenius-Perron dimension $pqr$ is
group-theoretical, where $p<q<r$ are distinct prime numbers.
\end{theorem}

Based on this theorem, Etingof et al then completed the classification of semisimple Hopf algebras of dimension $pqr$ and $pq^2$. They first proved the following lemma \cite[Lemma 9.5]{Etingof2}.

\begin{lemma}
Let $H$ be a group-theoretical semisimple Hopf algebra of square-free dimension.
Then $H$ fits into a split Abelian extension of the form $H(G,K,L, 1, 1, 1)$.
\end{lemma}

The definition and basic results on extensions of Hopf algebras can be easily found in the literature, e. g. \cite{Natale2}, \cite{Masuoka0}, \cite{Natale5}.

Let $p<q<r$ be prime numbers and $H$ be a semisimple Hopf algebra of dimension $pqr$. By the lemma above and Theorem \ref{thm9}, we have the following corollary \cite[Corollary 9.4]{Etingof2}.

\begin{corollary}
There exists a finite group $G$ of order $pqr$ and an exact factorization $G = KL$ of $G$
into a product of subgroups, such that $H$ is the split Abelian extension $H(G,K,L, 1, 1, 1) =k[K]\ltimes {\rm Fun}(L)$ associated to this factorization.
\end{corollary}

In \cite[Theorem 4.6]{Natale1}, Natale classified semisimple Hopf algebras of dimension $pqr$ which fit into Abelian extensions. Therefore, the corollary above gives a complete classification of semisimple Hopf algebras of dimension $pqr$.

Let $p,q$ be distinct prime numbers and $H$ be a semisimple Hopf algebra of dimension $pq^2$. By Theorem \ref{thm9}(1), the category Rep$H$ of finite-dimensional representations of $H$ is solvable. By Definition \ref{def1}, Rep$H$ is either an extension or an equivariantization of a fusion category of smaller dimension.  Etingof et al then proved that Rep$H$ is group-theoretical by considering these two possibilities \cite[Proposition 9.6]{Etingof2}. Consequently, they got the classification of semisimple Hopf algebras of dimension $pq^2$ as follows.

\begin{corollary}
A semisimple Hopf algebra of dimension $pq^2$ is either a Kac algebra, or a twisted
group algebra (by a twist corresponding to the subgroup $(\mathbb{Z}/q\mathbb{Z})^2$), or the dual of a twisted group
algebra.
\end{corollary}

\begin{remark}
Jordan and Larson \cite{Jordan} also proved, by different methods, that any semisimple Hopf algebra of dimension $pq^2$ is group-theoretical.

In a series of papers \cite{Natale1}, \cite{Natale0}, \cite{Natale3}, Natale studied the classification of semisimple Hopf algebras of dimension $pq^2$. In particular, she \cite{Natale3} completed the classification of semisimple Hopf algebras $H$ of dimension $pq^2$ such that $H$ and $H^*$ are both of Frobenius type. Therefore, Theorem \ref{thm8}, Theorem \ref{thm9} and Natale's results can also give the classification of semisimple Hopf algebras of dimension $pq^2$.
\end{remark}

Besides these applications, part (1) of Theorem \ref{thm9}
also provides a powerful method in classifying other semisimple Hopf
algebras whose dimensions consist of two prime divisors. For example,
let $H$ be a semisimple Hopf algebra of dimension $p^2q^2$, where
$p, q$ are distinct prime numbers with $p^4 < q$. Part (1) of
Theorem \ref{thm9} shows that the dimension of a simple $H$-module
can only be $1, p, p^2$ or $q$. It follows that we have an equation
$$p^2q^2 = |G(H^*)| + ap^2 + bp^4 + cq^2,$$
where $a, b, c$ is the number of non-isomorphic simple $H$-modules
of dimension $p, p^2$ and $q$, respectively. By analyzing the order
of $G(H^*)$ and standard subalgebras of $R(H)$, we can determine the
possible quotient Hopf algebras of $H$, and then obtain the classification of
$H$. See \cite{dong} for details.

\section{Further discussions}
To conclude this paper, we would like to discuss three questions which are tightly connected to Kaplansky's sixth conjecture.

As we have seen in the previous section, fusion category theory is a powerful tool
in the work toward solving Kaplansky's sixth conjecture. The Morita
equivalence method seems especially effective in this direction. This
is usually accomplished by analyzing the Drinfeld center of a fusion
category and then studying its Tannakian subcategories.

The following question is the second question in \cite{Etingof2}.  An negative answer to this question will solve Kaplansky's sixth conjecture, in view of Theorem \ref{thm8}.
\medbreak
\textbf{Question 1.}\, Does there exist a weakly integral fusion category which is not weakly group-theoretical?
\medbreak

Although the theory of Hopf algebras has developed for about $70$
years, we know little about the interrelations between Hopf algebras
and their duals. Let $H$ be a semisimple Hopf algebra over $k$, Rep$H$ and Rep$H^*$
be the category of finite-dimensional representations of $H$ and
$H^*$, respectively. The knowledge of interrelations
between Rep$H$ and Rep$H^*$ can greatly help us in solving
Kaplansky's sixth conjecture.

\medbreak
\textbf{Question 2.} For any semisimple Hopf algebra $H$, if $H$
satisfies Kaplansky's sixth conjecture, does $H^*$ satisfy Kaplansky's sixth conjecture?
\medbreak

If $H^*$ also satisfies Kaplansky's sixth conjecture then we can get closer to solving the conjecture. Since, by Theorem \ref{thm5}, $D(H)$ satisfies Kaplansky's
conjecture. Hence, by the assumption, $D(H)^*$ also satisfies Kaplansky's sixth conjecture. Therefore, the
dimension of every simple $H$-module divides the square of the dimension of $H$,
since $D(H)^*=H^{op}\otimes H^*$ as an algebra.

\medbreak
Let $H$ be a semisimple Hopf algebra and $\sigma: H\times H\to k$ be a normalized $2$-cocycle that is convolution-invertible, that is,
$$\sigma(x_1,y_1)\sigma(x_2y_2,z)=\sigma(y_1,z_1)\sigma(x,y_2z_2)\,\mbox{\,and\,}\sigma(1,1)=1,$$
where $x,y,z\in H$.

Let $H_{\sigma}=H$ as a coalgebra, but with the multiplication $._{\sigma}$ twisted by $\sigma$:
$$x._{\sigma}y=\sigma(x_1,y_1)x_2y_2\sigma^{-1}(x_3,y_3).$$
Then $H_{\sigma}$ is again a semisimple Hopf algebra. Moreover, $(H_{\sigma})_{\sigma^{-1}}=H$. We call $H_{\sigma}$ the twisting of $H$.

We do not know whether the class of finite-dimensional semisimple Hopf algebras is closed under twisting, in the positive characteristic setting. The reader is directed to \cite[Corollary 3.6 and Remark 3.9]{Aljadeff} for reference.

The above procedure is the dual version of twisting of coproduct which was introduced by Drinfeld \cite{Dri1986}. The reader can find a detailed exposition about these two twistings in \cite{Doi}.

\medbreak
\textbf{Question 3.} For any semisimple Hopf algebra $H$, if $H$
satisfies Kaplansky's sixth conjecture, does $H_{\sigma}$ satisfy Kaplansky's sixth conjecture?
\medbreak
If $H_{\sigma}$ also satisfies Kaplansky's sixth conjecture then the
dimension of every simple $H$-module divides the square of the dimension of $H$, because the Drinfeld double $D(H)$ is the twisting of $H^{*cop}\otimes H$ \cite{Doi}.

\medbreak
\textbf{{Acknowledgements:}}\quad The authors appreciate the suggestions and comments of the referee, which greatly improved this article. The authors are grateful to Huixiang Chen, Susan Montgomery and Shenglin Zhu for valuable discussions. The authors are grateful to Henry Tucker, who helped them with English language.


\begin{thebibliography}{*}
\bibitem{Aljadeff}A. Aljadeff, P. Etingof, S. Gelaki, D. Nikshych, \textit{On twisting of finite-dimensional Hopf algebras}, J. Algebra \textbf{256} (2002), pp. 484--501.

\bibitem{Andr}N. Andruskiewitsch, \textit{About finite dimensional Hopf algebras}, Contemp.
Math. \textbf{294} (2002), pp. 1--57.

\bibitem{BiNa2011}J. Bichon, S. Natale, \textit{Hopf algebra deformations of binary polyhedral groups}, Transform. Groups, \textbf{16} (2011), pp. 339--374.


\bibitem{Burciu}S. Burciu, \textit{Representations of degree three for semisimple
Hopf algebras}, J. Pure Appl. Algebra \textbf{194} (2004), pp. 303--312.

\bibitem{Burciu2}S. Burciu, \textit{On the classification of semisimple Hopf algebras: structure and applications}, Noncommutative structures in mathematics and physics, S. Caenepeel, J. Fuchs, A. Stolin (eds.), pp. 29--45, Proc. Royal Flemish Academy of Belgium, Brussels, 2009.

\bibitem{Curtis}C. W. Curtis, I. Reiner, \textit{Methods of representation theory with
applications to finite groups and orders}, Vol. I, Wiley Interscience, New York, 1981.


\bibitem{dong}J. Dong, \textit{Structure of semisimple Hopf algebras of dimension $p^2q^2$}, Comm. Algebra \textbf{40} (2012), pp. 795--806.


\bibitem{dong1}J. Dong, \textit{Structure theorems for semisimple Hopf algebras of dimension $pq^3$},  Comm. Algebra \textbf{40} (2012), pp. 4673--4678.

\bibitem{dodai}J. Dong, L. Dai, \textit{Further results on semisimple Hopf algebras
of dimension $p^2q^2$},  Rev. Un. Mat. Argentina \textbf{53} (2012), pp. 97--112.

\bibitem{dong2}J. Dong, S. Wang, \textit{On semisimple Hopf algebras of dimension $2q^3$}, J. Algebra \textbf{375} (2013), pp. 97--108.
\bibitem{dong2013integral}
J.~Dong, L.~Dai, \textit{On integral fusion categories with low-dimensional simple
  objects}, Comm. Algebra \textbf{42} (11), pp. 4955--4961.

\bibitem{dong2012frobenius}
J.~Dong, S.~Natale, L.~Vendramin, \textit{{F}robenius property for fusion categories of
  small integral dimension}, preprint arXiv:1209.1726, to appear in J.
  Algebra Appl..

\bibitem{Doi}Y. Doi, M. Takeuchi, \textit{Multiplication alteration by two-cocycles--the quantum version},
Comm. Algebra \textbf{22} (1994), pp. 5715--5732.

\bibitem{Dri1986}V. Drinfeld, \textit{Quantum groups}, Proceedings of the ICM, Berkeley, Amer. Math. Society, 1986.

\bibitem{Etingof}P. Etingof, S. Gelaki, \textit{Some properties of finite dimensional
semisimple Hopf algebras}, Math. Res. Lett. \textbf{5} (1998), pp. 191--197.

\bibitem{EtingofG}P. Etingof, S. Gelaki, \textit{Semisimple Hopf algebras of dimension $pq$ are trivial}, J. Algebra \textbf{210} (1998), pp. 664--669.


\bibitem{Etingof1}P. Etingof, D. Nikshych, V. Ostrik, \textit{On fusion categories}, Ann. Math. \textbf{162} (2005), pp. 581--642.

\bibitem{Etingof2}P. Etingof, D. Nikshych, V. Ostrik, \textit{Weakly group-theoretical and solvable fusion categories}, Adv. Math. \textbf{226} (2011), pp. 176--505.


\bibitem{GeWe}S. Gelaki, S. Westreich, \textit{On semisimple Hopf algebras of dimension
$pq$}, Proc. Amer. Math. Soc. \textbf{128} (2000), pp. 39--47. (Corrigendum: Proc. Amer. Math. Soc. \textbf{128} (2000), pp. 2829--2831)


\bibitem{Isaacs}M. Isaacs,\textit{ Character theory of finite groups}, Academic Press, New
York, 1976.

\bibitem{Jordan}D. Jordan, E. Larson, \textit{On the classification of certain fusion categories}, J. Noncommut. Geom. \textbf{3} (2009), pp. 481--499.

\bibitem{Kac}G. Kac, \textit{Certain arithmetic properties of ring groups}, Funct. Anal. Appl. \textbf{6} (1972), pp. 158--160.
\bibitem{KaPa1996}G. Kac, V. Paljutkin, \textit{Finite ring groups}, Trans. Moscow Math. Soc. \textbf{5} (1966), pp. 251--294.

\bibitem{Kaplansky}I. Kaplansky, \textit{Bialgebras}, University of Chicago Press, Chicago, 1975.

\bibitem{Kashina}Y. Kashina, Y. Sommerh\"{a}user, Y. Zhu, \textit{Self-dual modules of semisimple
Hopf algebras}, J. Algebra, \textbf{257} (2002), pp. 88--96.

\bibitem{Kashina2}Y. Kashina, Y. Sommerh\"{a}user, Y. Zhu, \textit{On higher Frobenius-Schur indicators}, Mem. Amer. Math. Soc. \textbf{181} (2006), 65 pp.

\bibitem{Kashina3}Y. Kashina, \textit{Some results on Hopf algebras of Frobenius type}, Math. Appl.
Dordr. \textbf{555} (2003), pp. 85--104.

\bibitem{Larson}R. G. Larson, D. E. Radford, \textit{Finite dimensional cosemisimple Hopf
algebras in characteristic $0$ are semisimple}, J. Algebra \textbf{117} (1988),  pp. 267--289.

\bibitem{Larson2}R. G. Larson, D. E. Radford, \textit{Semisimple
cosemisimple Hopf algebras}, Amer. J. Math. \textbf{109} (1987),  pp. 187--195.


\bibitem{Linchenko}V. Linchenko, S. Montgomery, \textit{A Frobenius-Schur theorem for Hopf
algebras}, Algebras Represent. Theory \textbf{3} (2000),  pp. 347--355.

\bibitem{Lorenz}M. Lorenz, \textit{On the degree of irreducible representations of Hopf algebras}, arxiv:math/0003008.

\bibitem{Lorenz2}M. Lorenz, \textit{On the class equation for Hopf algebras},
Proc. Amer. Math. Soc. \textbf{126} (1998),  pp. 1841--2844.

\bibitem{majid}S. Majid, \textit{Doubles of quasitriangular Hopf algebras}, Comm. Algebra \textbf{19} (1991),  pp. 3061--3073.
\bibitem{Masuoka0}A. Masuoka, \textit{Extensions of Hopf algebras},  Trabajos de Matem¨¢tica 41/99, FaMAF.

\bibitem{Masuoka1}A. Masuoka, \textit{The $p^n$ theorem for semisimple Hopf
algebras}, Proc. Amer. Math. Soc. \textbf{124} (1996),  pp. 735--737.

\bibitem{Masuoka2}A. Masuoka, \textit{Self-dual Hopf algebras of dimension $p^3$ obtained by
extension}, J. Algebra \textbf{178} (1995),  pp. 791--806.

\bibitem{Montgomery}S. Montgomery, \textit{Hopf algebras and their actions on rings}, CBMS Reg.
Conf. Ser. Math. 82, Amer. Math. Soc., Providence, 1993.

\bibitem{Montgomery2}S. Montgomery, S. Witherspoon, \textit{Irreducible representations of
crossed products}, J. Pure Appl. Algebra \textbf{129} (1998),  pp. 315--326.


\bibitem{Natale1}S. Natale, \textit{On semisimple Hopf algebras of dimension $pq^2$}, J. Algebra \textbf{221} (1999),  pp. 242--278.
\bibitem{Natale0}S. Natale, \textit{On semisimple Hopf algebras of dimension $pq^2$, II}, Algebras Represent. Theory \textbf{4} (2001),  pp. 277--291.
\bibitem{Natale2}S. Natale, \textit{On group theoretical Hopf algebras and exact factorizations of finite groups,} J. Algebra \textbf{270} (2003),  pp. 199--211.

\bibitem{Natale3}S. Natale, \textit{On semisimple Hopf algebras of dimension $pq^r$}, Algebras Represent. Theory \textbf{7} (2004),  pp. 173--188.

\bibitem{Natale4}S. Natale, \textit{Semisolvability of semisimple Hopf algebras of low dimension}, Mem. Amer. Math. Soc. \textbf{186} (2007),  123 pp.

\bibitem{Natale5}S. Natale, \textit{Semisimple Hopf algebras and their representations}, Publ. Mat. Urug. \textbf{12} (2011),  pp. 123--167.

\bibitem{Nichols}W. D. Nichols, M. B. Richmond, \textit{The Grothendieck group of a Hopf
algebra}, J. Pure Appl. Algebra \textbf{106} (1996),  pp. 297--306.

\bibitem{Nichols2}W. D. Nichols, M. B. Zoeller, \textit{A Hopf algebra freeness theorem}, Amer.
J. Math. \textbf{111} (1989),  pp. 381--385.

\bibitem{Sch}H.-J. Schneider, \textit{Some properties of factorizable Hopf
algebras}, Proc. Amer. Math. Soc. \textbf{129} (2001),  pp. 1891--1898

\bibitem{Sch2}H.-J. Schneider,
\textit{Normal basis and transitivity of crossed products for Hopf algebras},
J. Algebra \textbf{152} (1992),  pp. 289--312.

\bibitem{Sommer}Y. Sommerh\"{a}user, \textit{On Kaplansky's conjectures}, Interactions between ring theory and representations of algebras, F. v. Oystaeyen, M. Saorin (eds.),  pp. 393--412, Lecture Notes in Pure and Appl. Math. 210, Dekker, New York, 2000.

\bibitem{Somm}Y. Sommerh\"{a}user, \textit{Yetter-Drinfel¡¯d Hopf algebras over groups of
prime order}, Lecture Notes in Math. 1789, Springer-Verlag, 2002.

\bibitem{Sweedler}M. E. Sweedler, \textit{Integrals for Hopf algebras}, Ann. Math. \textbf{89} (1969), pp. 323--335.

\bibitem{Sweedler2}M. E. Sweedler, \textit{Hopf Algebras}, Benjamin, New York, 1969.
\bibitem{Tsang}Y. Tsang, Y. Zhu, \textit{On the Drinfeld double of a Hopf
algebra}, preprint, 1998.


\bibitem{YaZh}M. Yan, Y. Zhu, \textit{Stabilizer for Hopf algebra actions}, Comm. Algebra \textbf{26} (1998),  pp. 3885--3898.

\bibitem{ZhuS}S. Zhu, \textit{On finite dimensional semisimple Hopf algebras}, Comm. Algebra \textbf{21} (1993),  pp. 3871--3885.

\bibitem{Zhu}Y. Zhu, \textit{Hopf algebras of prime dimension}, Internat. Math. Res. Notices
\textbf{1 }(1994),  pp. 53--59.

\bibitem{Zhu2}Y. Zhu, \textit{A commuting pair in Hopf algebras}, Proc. Amer. Math. Soc.
\textbf{125} (1997),  pp. 2847--2851.


\bibitem{Zhu3}Y. Zhu, \textit{The dimension of irreducible modules for transitive module algebras}, Comm. Algebra \textbf{29} (2001),  pp. 2877--2886.
\end{thebibliography}
\end{document}